\documentclass[12pt]{article}
\usepackage{amssymb}
\usepackage{amsfonts}
\usepackage{amsmath}
\usepackage{amsthm}

\newtheorem{theorem}{Theorem}[section]
\newtheorem{corollary}[theorem]{Corollary}
\newtheorem{lemma}[theorem]{Lemma}
\newtheorem{proposition}[theorem]{Proposition}
\theoremstyle{definition}

\newtheorem{conjecture}[theorem]{Conjecture}

\newtheorem{definition}[theorem]{Definition}
\newtheorem{example}[theorem]{Example}

\newtheorem{remark}[theorem]{Remark}

\input{tcilatex}
\begin{document}

\title{A unified approach to the plus-construction, Bousfield localization,
Moore spaces and zero-in-the-spectrum examples}
\author{Shengkui Ye}
\maketitle

\begin{abstract}
We introduce a construction adding low-dimensional cells to a space that
satisfies certain low-dimensional conditions; it preserves high-dimensional
homology with appropriate coefficients. This includes as special cases
Quillen's plus construction, Bousfield's integral homology localization, the
existence of Moore spaces $M(G,1)$ and Bousfield and Kan's partial $k$%
-completion of spaces. We also use it to generalize counterexamples to the
zero-in-the-spectrum conjecture found by Farber and Weinberger, and by
Higson, Roe and Schick.
\end{abstract}

\section{Introduction}

The aim of this article is to give a unified treatment of Quillen's
plus-construction, Rodr\'{\i}guez and Scevenels' work on Bousfield's
integral localization, Varadarajan's theorem on the existence of Moore
spaces, the partial $k$-completion of Bousfield and Kan, and counterexamples
to the zero-in-the-spectrum conjecture by Farber and Weinberger, and Higson,
Roe and Schick. We introduce a method for adding low-dimensional cells to a
space satisfying certain low-dimensional conditions and preserving
high-dimensional homology with appropriate coefficients. First, we briefly
review these existing works.

Let $X$ be a CW complex with fundamental group $G$ and $P$ a perfect normal
subgroup of $G$, \textsl{i.e.}, $P=[P,P].$ Quillen \cite{qu} shows that
there exists a CW complex $X_{P}^{+}$, whose fundamental group is $G/P,$ and
an inclusion $f:X\rightarrow X_{P}^{+}$ such that 
\begin{equation*}
H_{n}(X;f_{\ast }M)\cong H_{n}(X_{P}^{+};M)
\end{equation*}%
for any integer $n$ and local coefficient system $M$ over $X_{P}^{+}.$ Here $%
X_{P}^{+}$ is called the plus-construction of $X$ with respect to $P$ and is
unique up to homotopy equivalence. The plus-construction can be used to
define higher algebraic $K$-theory, as follows. Let $R$ be a unital
associative ring with $n$-th general linear group $\mathrm{GL}_{n}(R)$ and $%
E_{n}(R)$ its subgroup generated by elementary matrices. Let $\mathrm{GL}%
(R)=\cup _{n}\mathrm{GL}_{n}(R)$. Then $E(R)=\cup _{n}E_{n}(R)$ is the
maximal perfect normal subgroup of $\mathrm{GL}(R)$. For the classifying
space $B\mathrm{GL}(R),$ let $B\mathrm{GL}(R)^{+}$ denote its
plus-construction with respect to $E(R)$. The algebraic $K$-groups are
defined as $K_{i}(R)=\pi _{i}(B\mathrm{GL}(R)^{+})$ ($i\geq 1$).

While the plus-construction preserves all ordinary homology with
coefficients, for each generalized homology theory there is a Bousfield
localization preserving it. More precisely, write \textsc{Ho} for the
pointed homotopy category of CW complexes. Then Bousfield \cite{Bs} shows
that each generalized homology theory $h_{\ast }$ determines an $h_{\ast }$%
-localization functor $E:$ \textsc{Ho }$\rightarrow $ \textsc{Ho} and a
natural transformation $\eta :\mathrm{Id}\rightarrow E$. This localization
is characterized by the universal property that $\eta _{X}:X\rightarrow EX$
is the terminal $h_{\ast }$-homology equivalence going out of $E$, \textsl{%
i.e.,}

\begin{enumerate}
\item[(i)] $\eta _{X}:X\rightarrow EX$ induces $h_{\ast }(X)\cong h_{\ast
}(EX)$, and

\item[(ii)] for any map $f:X\rightarrow Y\in $ \textsc{Ho} inducing $h_{\ast
}(X)\cong h_{\ast }(Y)$ there is a unique map $r:Y\rightarrow EX\in $ 
\textsc{Ho} with $rf=\eta _{X}.$
\end{enumerate}

For ordinary homology theory $H\mathbb{Z}$ with $\mathbb{Z}$ as
coefficients, Bousfield's $H\mathbb{Z}$-localization $X_{H\mathbb{Z}}$ of a
space $X$ is homotopy equivalent to localization of $X$ with respect to a
map of classifying spaces $Bf:BF_{1}\rightarrow BF_{2}$ induced by a certain
homomorphism $f:F_{1}\rightarrow F_{2}$ of free groups (\textsl{cf}. \cite%
{8,9}). This implies that a space $X$ is $H\mathbb{Z}$-local if and only if
the induced map $Bf^{\ast }:\mathrm{map}(BF_{2},X)\rightarrow $\textrm{map}$%
(BF_{1},X)$ is a weak homotopy equivalence. Rodr\'{\i}guez and Scevenels 
\cite{RS} show that there is a topological construction that, while leaving
the integral homology of a space unchanged, kills the intersection of the
transfinite lower central series of its fundamental group. Moreover, this is
the maximal subgroup that can be factored out of the fundamental group
without changing the integral homology of a space. For more information on $H%
\mathbb{Z}$-localization and homology equivalence with other coefficients,
see \cite{1,3,4,me,mp} and references therein.

The plus-construction shares some common features with the construction of
Moore spaces in \cite{ann}. Given an integer $n\geq 1$ and a group $G$
(abelian if $n\geq 2$), a Moore space $M(G,n)$ is a space $X$ such that $\pi
_{j}(X)=0$ for $j<n$, $\pi _{n}(X)=G$ and $H_{i}(X;\mathbb{Z})=0$ for $i>n.$
For $n\geq 2,$ such a space always exists. For $n=1,$ Varadarajan \cite{ann}
proves that there exists a Moore space $M(G,1)$ if and only if $H_{2}(G;%
\mathbb{Z})=0.$

Let $k$ be the constant ring $\mathbb{Z}/p$ (prime $p$) or $k\subseteq 
\mathbb{Q}$ a subring of the rationals. For a space $X,$ let $P_{\pi
_{1}(X)} $ be the largest subgroup of $\pi _{1}(X)$ for which $H_{1}(P_{\pi
_{1}(X)};k)=0.$ Bousfield and Kan \cite{Bk} show that there exists a space $%
C^{k}(X),$ which is called the partial $k$-completion of $X,$ and a map $%
\phi :X\rightarrow C^{k}(X)$ such that the fundamental group $\pi
_{1}(C^{k}(X))=\pi _{1}(X)/P_{\pi _{1}(X)}$ and for any integer $q\geq 0\ $%
the map $\phi $ induces an isomorphism $H_{q}(X;k[\pi /P])\cong
H_{q}(Y;k[\pi /P]).$ (Here, $k[\pi /P]$ is the group ring over $k$ of $\pi
/P.$)

The zero-in-the-spectrum conjecture goes back to Gromov, who asked whether
for a closed, aspherical, connected and oriented Riemannian manifold $M$
there always exists some $p\geq 0,$ such that zero belongs to the spectrum
of the Laplace-Beltrami operator $\Delta _{p}$ acting on the square
integrable $p$-forms on the universal covering $\tilde{M}$ of $M.$ Farber
and Weinberger \cite{FW} show that the conjecture is not true if the
condition that $M$ is aspherical is dropped. More generally, Higson, Roe and
Schick \cite{hig} show that for a finitely presented group $G$ satisfying $%
H_{0}(G;C_{r}^{\ast }(G))=H_{1}(G;C_{r}^{\ast }(G))=H_{2}(G;C_{r}^{\ast
}(G))=0,$ there always exists a finite CW complex $Y$ with $\pi _{1}(Y)=G$
such that $Y$ is a counterexample to the conjecture if $M$ is not required
to be aspherical.

In this note, a more general construction is provided to preserve homology
theory. For this, we have to introduce the notion of a $G$-dense ring (for
details, see Definition \ref{dense}). Examples of $G$-dense rings include
the real reduced group $C^{\ast }$-algebra $C_{\mathbb{R}}^{\ast }(G),$ the
real group von Neumann algebra $\mathcal{N}_{\mathbb{R}}G,$ the real Banach
algebra $l_{\mathbb{R}}^{1}(G),$ the constant rings $k=\mathbb{Z}/p$ (prime $%
p$) and $k\subseteq \mathbb{Q}$ a subring of the rationals, the group ring $%
k[G],$ and so on.

\textit{Notation.} Let $\pi $ and $G$ be two groups. Suppose $R$ is a $%
\mathbb{Z}[G]$-module and $BG,B\pi $ are the classifying spaces. For a group
homomorphism $\alpha :\pi \rightarrow G,$ we will denote by $H_{1}(G,\pi ;R)$
the relative homology group $H_{1}(BG,B\pi ;R)$ with coefficients $R.$ All
spaces are assumed to be connected.

\begin{theorem}
\textit{\ \label{th1}}Assume that $G$ is a group and $(R,\phi )$ is a $G$%
-dense ring. Let $X$ be a CW complex with fundamental group $\pi =\pi
_{1}(X) $. Assume $\alpha :\pi \rightarrow G$ is a group homomorphism such
that 
\begin{eqnarray*}
\alpha _{\ast } &:&H_{1}(\pi ;R)\rightarrow H_{1}(G;R)\text{ is injective,
and } \\
\alpha _{\ast } &:&H_{2}(\pi ;R)\rightarrow H_{2}(G;R)\text{ is surjective.}
\end{eqnarray*}%
Suppose either that $R$ is a principal ideal domain or that the relative
homology group $H_{1}(G,\pi ;R)$ is a stably free $R$-module. Then there
exists a CW complex $Y$ with the following properties:

\begin{enumerate}
\item[(i)] $Y$ is obtained from $X$ by adding 1-cells, 2-cells and 3-cells,
such that

\item[(ii)] $\pi _{1}(Y)=G$ and the inclusion map $g:X\rightarrow Y$ induces
the same fundamental group homomorphism as $\alpha ,$ and

\item[(iii)] for any $q\geq 2$ the map $g$ induces an isomorphism%
\begin{equation}
g_{\ast }:H_{q}(X;R)\overset{\cong }{\rightarrow }H_{q}(Y;R).  \tag{1}
\end{equation}
\end{enumerate}
\end{theorem}

\bigskip

Theorem \ref{th1} has many important consequences, including the
following.\medskip

\begin{enumerate}
\item[(1)] \noindent $\alpha $ \textbf{surjective}.
\end{enumerate}

\begin{itemize}
\item When $R=\mathbb{Z}$ and $\ker \alpha $ is perfect, Proposition \ref%
{prop1} shows that Quillen's plus-construction is a special case of this
theorem.

\item When $R=\mathbb{Z}$, we obtain the result of Rodr\'{\i}guez and
Scevenels \cite{RS} on Bousfield integral localization (\textsl{cf}.
Corollary \ref{cor}).

\item When $k=\mathbb{Z/}p$ or $k\subseteq \mathbb{Q}$ a subring of the
rationals and $R=k[G]$, the theorem yields the partial $k$-completion of
Bousfield and Kan \cite{Bk} (see Corollary \ref{cmp}).
\end{itemize}

\begin{enumerate}
\item[(2)] \noindent $\pi =1.$
\end{enumerate}

\begin{itemize}
\item When $R=\mathbb{Z}$, we obtain in Corollary \ref{cor2} the existence
of the Moore space $M(G,1),$ which was first proved by Varadarajan in \cite%
{ann}.

\item When $R=C_{\mathbb{R}}^{\ast }(G),$ the theorem yields the results
obtained by Farber-Weinberger \cite{FW} and Higson-Roe-Schick \cite{hig} on
the zero-in-the-spectrum conjecture (actually our result is more general,
see Corollary \ref{zero}).
\end{itemize}

In Section 2, we give the definition and some examples of $G$-dense rings.
The main result Theorem \ref{th1} is proved in Section 3. Several
applications are presented in the last section.

\section{$G$-dense rings\label{s1}}

Let $G$ be a group. In this section, we introduce a kind of rings which
includes as special cases the real reduced group $C^{\ast }$-algebra $C_{%
\mathbb{R}}^{\ast }(G),$ the real group von Neumann algebra $\mathcal{N}_{%
\mathbb{R}}G,$ the real Banach algebra $l_{\mathbb{R}}^{1}(G),$ $k=\mathbb{Z/%
}p$ for some prime $p$ or $k\subseteq \mathbb{Q}$ a subring of the
rationals, and the group ring $k[G]$.

\begin{definition}
\label{dense}A $G$\emph{-dense ring} $(R,\phi )$ is a unital ring $R$
together with a ring homomorphism $\phi :\mathbb{Z}[G]\rightarrow R$ such
that, when $R$ is regarded as a left $\mathbb{Z}[G]$-module via $\phi ,$
then, for any right $\mathbb{Z}[G]$-module $M,$ free right $R$-module $F$
and $R$-module surjection $f:M\bigotimes_{\mathbb{Z}[G]}R\twoheadrightarrow
F,$ the module $F$ has an $R$-basis in $f(M\bigotimes 1).$
\end{definition}

When $\phi $ is obvious, it is omitted from the notation. Some examples of $%
G $-dense rings are as follows. Recall that for a group $G$, the space 
\begin{equation*}
l^{2}(G)=\{f:G\rightarrow \mathbb{C}\mid \sum_{g\in G}|f(g)|^{2}<+\infty \}
\end{equation*}
is a Hilbert space with inner product $\langle f_{1},f_{2}\rangle
=\sum_{x\in G}f_{1}(x)\overline{f_{2}(x)}$. Let $B(l^{2}(G))$ be the set of
all bounded linear operators of the Hilbert space $l^{2}(G).$ By definition,
the real reduced group $C^{\ast }$-algebra $C_{\mathbb{R}}^{\ast }(G)$ is
the completion of $\mathbb{R}[G]$ in $B(l^{2}(G))$ with respect to the
operator norm, while the real group von Neumann algebra $\mathcal{N}_{%
\mathbb{R}}G$ is the completion of $\mathbb{R}[G]$ in $B(l^{2}(G))$ with
respect to the weak operator norm. The real Banach algebra $l_{\mathbb{R}%
}^{1}(G)$ is the completion of the group ring $\mathbb{R}[G]$ with respect
to the $l^{1}$-norm.

\begin{lemma}
\label{lm2}The set of $G$-dense rings contains the real reduced group $%
C^{\ast }$-algebra $C_{\mathbb{R}}^{\ast }(G),$ the real group von Neumann
algebra $\mathcal{N}_{\mathbb{R}}G,$ the real Banach algebra $l_{\mathbb{R}%
}^{1}(G),$ the constant rings $k=\mathbb{Z/}p$ for any prime $p$, $%
k\subseteq \mathbb{Q}$ a subring of the rationals, and the group rings $%
k[G]. $
\end{lemma}

\begin{proof}
We prove the lemma case by case. Let $M$ be a right $\mathbb{Z}[G]$-module, $%
F$ a free right $R$-module and $f:M\bigotimes_{\mathbb{Z}[G]}R%
\twoheadrightarrow F$ a surjection of $R$-modules. Choose a basis $%
(b_{i})_{i\in S}$ of $F$ for some index set $S$. Since $f$ is surjective and 
$R$-linear, we can assume%
\begin{equation*}
b_{i}=\sum f(x_{ik}\tbigotimes 1)a_{ik}
\end{equation*}%
for some $x_{ik}\in M$ and $a_{ik}\in R.$

\begin{enumerate}
\item[(i)] \noindent $k=Z/p$\textbf{\ for some prime }$p,$\textbf{\ and }$%
R=k $\textbf{\ or} $k[G].$
\end{enumerate}

The ring homomorphism $\phi :\mathbb{Z}[G]\rightarrow R$ is induced from the
natural map $\mathbb{Z\rightarrow Z/}p.$ There is a surjection $\beta :%
\mathbb{Z}[G]\rightarrow R.$ Choose some $\tilde{a}_{ik}\in \beta
^{-1}(a_{ik}).$ Then we have 
\begin{equation*}
b_{i}=\sum f(x_{ik}\cdot \tilde{a}_{ik}\tbigotimes 1),
\end{equation*}%
which is in the image of $f(M\bigotimes 1).$

\begin{enumerate}
\item[(ii)] \noindent $k\subseteq Q$ \textbf{a subring of the rationals, and 
}$R=k$\textbf{\ or} $k[G].$
\end{enumerate}

The ring homomorphism $\phi :\mathbb{Z}[G]\rightarrow R$ is induced from the
natural map $\mathbb{Z}\rightarrow k.$ There is an inclusion $\beta
:R\hookrightarrow \mathbb{Q}[G].$ Then there exists an integer $n_{i},$
which is invertible in $R,$ such that $n_{i}a_{ik}\in \mathbb{Z}[G]$ and%
\begin{equation*}
n_{i}b_{i}=\sum f(x_{ik}\tbigotimes 1)n_{i}a_{ik}=\sum f(x_{ik}\cdot
n_{i}a_{ik}\tbigotimes 1),
\end{equation*}%
which is in the image of $f(M\bigotimes 1).$ Since $n_{i}$ is invertible, $%
(n_{i}b_{i})_{i\in S}$ is still a basis.

\begin{enumerate}
\item[(iii)] \noindent $R=C_{\mathbb{R}}^{\ast }(G),$ $l_{\mathbb{R}}^{1}(G)$%
\textbf{\ or }$\mathcal{N}_{\mathbb{R}}G.$
\end{enumerate}

The ring homomorphism $\phi :\mathbb{Z}[G]\rightarrow R$ is the natural
inclusion. The proof is similar to that of Proposition 4.4 in \cite{hig}. We
just briefly repeat here. Assume $R=C_{\mathbb{R}}^{\ast }(G)$, while the
other cases are similar. Choose $a_{ik}\in C_{\mathbb{R}}^{\ast }(G).$ Since 
$F$ is a free $C_{\mathbb{R}}^{\ast }(G)$ module, there is a natural product
topology on $F.$ As the set of all bases in $F$ is open and the module
multiplication operation%
\begin{equation*}
F\times C_{\mathbb{R}}^{\ast }(G)\rightarrow F
\end{equation*}%
is continuous, for each pair $(i,k)$ we can choose $a_{ik}^{\prime }\in 
\mathbb{Q}[G]$ sufficiently close to $a_{ik}$ such that the elements 
\begin{equation*}
b_{i}^{\prime }=\sum f(x_{ik})\tbigotimes a_{ik}^{\prime }
\end{equation*}%
form a new basis for $F.$ Since the tensor is over $\mathbb{Z}[G],$ a
similar argument as in the case of $k\subseteq \mathbb{Q}$, and $R=k$ or $%
k[G]$ shows that the image of $f(M\bigotimes 1)$ contains a basis.
\end{proof}

The following lemma provides more examples of $G$-dense rings.

\begin{lemma}
Let $G$ be a group and $N$ a normal subgroup of $G$ inducing the canonical
surjection $\psi :\mathbb{Z}[G]\rightarrow \mathbb{Z}[G/N].$ Assume that $%
(R,\phi )$ is a $G/N$-dense ring. Then $(R,\phi \circ \psi )$ is $G$-dense.
\end{lemma}

\begin{proof}
Suppose that, for a right $\mathbb{Z}[G]$-module $M$ and a free right $R$%
-module $F,$ there is a right $R$-module surjection $f:M\bigotimes_{\mathbb{Z%
}[G]}R\twoheadrightarrow F.$ Since 
\begin{equation*}
M\tbigotimes\nolimits_{\mathbb{Z}[G]}R\cong (M\tbigotimes\nolimits_{\mathbb{Z%
}[G]}\mathbb{Z[}G/N])\tbigotimes\nolimits_{\mathbb{Z}[G\mathbb{/}N]}R,
\end{equation*}%
the module $F$ has an $R$-basis in $f(M\bigotimes_{\mathbb{Z}[G]}\mathbb{Z[}%
G/N]\bigotimes_{\mathbb{Z[}G/N]}1).$ The quotient map $\psi $ gives 
\begin{equation*}
f(M\tbigotimes\nolimits_{\mathbb{Z}[G]}\mathbb{Z[}G/N]\tbigotimes\nolimits_{%
\mathbb{Z[}G/N]}1)
\end{equation*}%
as a subset of $f(M\bigotimes_{\mathbb{Z}[G]}1).$
\end{proof}

Recall from \cite{nw,cn1,cn2,le} that a ring $R$ is \textit{right Steinitz}
if $R$ has the property that any linearly independent subset of a free right 
$R$-module $F$ can be extended to a basis of $F$. It is also known that
right Steinitz rings are precisely the right perfect local rings (\textsl{cf.%
} \cite{nw}).

\begin{proposition}
Let $G$ be a group. Any right Steinitz ring $R$ with a ring homomorphism $%
\phi :\mathbb{Z}[G]\rightarrow R$ is $G$-dense.
\end{proposition}

\begin{proof}
Let $R$ be a right Steinitz ring, $M$ a right $\mathbb{Z}[G]$-module, $F$ a
free right $R$-module and $f:M\bigotimes_{\mathbb{Z}[G]}R\twoheadrightarrow
F $ an $R$-module surjection. By Theorem 1.1 in \cite{nw}, any generating
set of a free right $R$-module $F$ contains a basis of $F$. Since the set $%
f(M\bigotimes 1)$ is a generating set of $F,$ this shows that $R$ is $G$%
-dense.
\end{proof}

Before we give an example of a ring which is not $G$-dense, let's present a
matrix property of $G$-dense rings. For a matrix $A=(a_{ij})\in M_{n}(R),$
denote by $A_{k}$ the submatrix spanned by the first $k$ columns.

\begin{proposition}
\label{gd}Let $G$ be a group, $R$ a unital ring with invariant basis number (%
\textrm{IBN}, \textsl{cf.} \cite{mag}) and $n$ a positive integer. Assume
that $(R,\phi )$ is $G$-dense. Then for every integer $n\geq 1$, any matrix $%
A\in \mathrm{GL}_{n}(R)$ and any integer $k$ with $1\leq k\leq n,$ there
exists a matrix $B\in M_{k\times n}(\mathbb{Z}[G])$ such that $\phi
(B)A_{k}\in \mathrm{GL}_{k}(R).$
\end{proposition}

\begin{proof}
Let $n\geq 1$ be an integer and the standard basis of $R^{n}$ be denoted as $%
\{\mathbf{f}_{1},\mathbf{f}_{2},\ldots ,\mathbf{f}_{n}\}.$ Write the element 
$r_{1}\mathbf{f}_{1}+r_{2}\mathbf{f}_{2}+\cdots +r_{n}\mathbf{f}_{n}$ in $%
R^{n}$ as the vector $[r_{1},r_{2},\ldots ,r_{n}].$ Let $\{\mathbf{e}_{1},%
\mathbf{e}_{2},\ldots ,\mathbf{e}_{n}\}$ be the standard basis of $(\mathbb{Z%
}[G])^{n}.$ Given $A=(a_{ij})\in \mathrm{GL}_{n}(R)$, define the map 
\begin{equation*}
\alpha :(\mathbb{Z}[G])^{n}\tbigotimes\nolimits_{\mathbb{Z}[G]}R\rightarrow
R^{n}
\end{equation*}%
by $\alpha (\sum \mathbf{e}_{i}\bigotimes r_{i})=[r_{1},r_{2},\ldots
,r_{n}]A.$ For each integer $k$ such that $1\leq k\leq n,$ let $p$ be the
standard projection to the first $k$ components%
\begin{equation*}
p:R^{n}\rightarrow R^{k}.
\end{equation*}%
Since $A$ is invertible, we have a surjection%
\begin{equation*}
p\circ \alpha :(\mathbb{Z}[G])^{n}\tbigotimes\nolimits_{\mathbb{Z}%
[G]}R\rightarrow R^{k}.
\end{equation*}
According to the definition of $G$-dense rings and the assumption that $R$
has \textrm{IBN}, we can find elements $\mathbf{x}_{1},\mathbf{x}_{2},\ldots
,\mathbf{x}_{k}\in (\mathbb{Z}[G])^{n}$ such that $\{p\circ \alpha (\phi (%
\mathbf{x}_{i}))=\phi (\mathbf{x}_{i})A\mid i=1,2,\ldots ,k\}$ is a basis
for $R^{k}.$ Let $B=[\mathbf{x}_{1},\mathbf{x}_{2},\ldots ,\mathbf{x}%
_{k}]^{T}.$ By the definitions of $\alpha $ and $p,$ the matrix $\phi
(B)A_{k}$ is invertible$.$
\end{proof}

According to Proposition \ref{gd}, the following example shows that the ring
of Gauss integers $\mathbb{Z}[i]$ is not $G$-dense for the trivial group $G.$

\begin{example}
Let $\mathbb{Z[}i\mathbb{]}$ be the Gauss integers. Note that the matrix%
\begin{equation*}
\begin{bmatrix}
3 & 2-i \\ 
i+2 & 2%
\end{bmatrix}%
\end{equation*}%
lies in $\mathrm{SL}_{2}(\mathbb{Z[}i\mathbb{]})$. But we are not able to
find two integers $a,b$ such that $3a+(i+2)b$ is invertible in $\mathbb{Z[}i%
\mathbb{]}$, since the only units are $1,-1,i,-i.$ This shows $\mathbb{Z[}i%
\mathbb{]}$ is not $G$-dense for the trivial group $G$.
\end{example}

\section{The generalized plus-construction}

In this section, we will prove the main result, Theorem \ref{th1}. Theorem %
\ref{th1} shows that for certain homology theories, there is a construction
that preserves the higher homology groups.

In order to prove Theorem \ref{th1}, we use the following lemma, which is a
more general version of Hopf's exact sequence.

\begin{lemma}[Lemma 2.2 in \protect\cite{hig}]
\label{lem1}Let $G$ be a group and $V$ be a left $\mathbb{Z}[G]$-module. For
any CW complex $X$ with fundamental group $G$ and universal covering space $%
\tilde{X}$, there is an exact sequence%
\begin{equation*}
H_{2}(\tilde{X})\tbigotimes\nolimits_{\mathbb{Z}[G]}V\rightarrow
H_{2}(X;V)\rightarrow H_{2}(G;V)\rightarrow 0.
\end{equation*}
\end{lemma}

\begin{proof}[Proof of Theorem \protect\ref{th1}]
For the group homomorphism $\alpha :\pi _{1}(X)=\pi \rightarrow G,$ we
construct a CW complex $W$ such that $\pi _{1}(W)=G$ as follows. Let $S$ be
a set of normal generators of $\ker (\alpha ),$ \textsl{i.e.}, $\ker (\alpha
)$ is generated by elements of the form $gsg^{-1}$ for $s\in S$ and $g\in
\pi .$ For each element in $S$, attach a 2-cell $(D^{2},S^{1})$ to $X$ to
kill the corresponding element in $\pi $. Extend the presentation of $\pi /%
\mathrm{ker}(\alpha )$ by generators and relations to a presentation of $G$.
This can be obtained by adding generators and relations to the presentation
of $\pi /\mathrm{ker}(\alpha )$. When $f$ is surjective, we only need to add
relations. For each such generator (resp. relation), we continue to add a
1-cell (resp. 2-cell) to $X$ (see 5.1 in \cite{berr} for more details). Let $%
W$ denote the resulting space$.$

We consider the homology groups of the pair $(W,X).$ By Lemma \ref{lem1},
there is a commutative diagram%
\begin{equation*}
\begin{array}{cccc}
H_{2}(\tilde{X})\tbigotimes\limits_{\mathbb{Z}[G]}R & \rightarrow & H_{2}(%
\tilde{W})\tbigotimes\limits_{\mathbb{Z}[G]}R &  \\ 
\downarrow &  & \downarrow j_{4} &  \\ 
H_{2}(X;R) & \overset{j_{2}}{\longrightarrow } & H_{2}(W;R) & \overset{j_{1}}%
{\rightarrow }H_{2}(W,X;R)\rightarrow H_{1}(X;R)\rightarrow H_{1}(W;R) \\ 
\downarrow j_{3} &  & \downarrow j_{5} &  \\ 
H_{2}(\pi ;R) & \overset{\alpha _{\ast }}{\longrightarrow } & H_{2}(G;R) & 
\end{array}%
\end{equation*}%
where the middle horizontal chain is the long exact sequence of homology
groups for the pair $(W,X)$ and the two vertical lines are the exact
sequences as in Lemma \ref{lem1}. Notice that 
\begin{equation*}
H_{1}(X;R)\cong H_{1}(\pi ;R)\rightarrow H_{1}(W;R)\cong H_{1}(G;R)
\end{equation*}%
is injective by assumption. This implies $j_{1}:H_{2}(W;R)\rightarrow
H_{2}(W,X;R)$ is surjective in the above diagram. For any element $a\in
H_{2}(W,X;R)$, choose its preimage $b\in H_{2}(W;R)$. Since $\alpha _{\ast
}:H_{2}(\pi ;R)\rightarrow H_{2}(G;R)$ is surjective by assumption, there
exists some element $c\in H_{2}(\pi ;R)$ such that $\alpha _{\ast
}(c)=j_{5}(b)$. Let $d\in H_{2}(X;R)$ be a preimage of $c,$ \textsl{i.e.} $%
j_{3}(d)=c$. By the commutativity of the diagram, $j_{5}(b-j_{2}(d))=0$.
Therefore, there exists an element $e\in H_{2}(\tilde{W})\tbigotimes_{%
\mathbb{Z}[G]}R$ such that $j_{4}(e)=b-j_{2}(d).$ It can be checked that $%
j_{1}\circ j_{4}(e)=a.$ Therefore, there is a surjection 
\begin{equation*}
j_{1}\circ j_{4}:H_{2}(\tilde{W})\tbigotimes\nolimits_{\mathbb{Z}%
[G]}R\rightarrow H_{2}(W,X;R)
\end{equation*}%
by this diagram chase.

We show that the relative homology group $H_{2}(W,X;R)$ can be taken to be a
free $R$-module. Let 
\begin{equation*}
0\rightarrow C_{2}(\tilde{W},\tilde{X};R)\overset{i}{\rightarrow }C_{1}(%
\tilde{W},\tilde{X};R)\overset{j}{\rightarrow }C_{0}(\tilde{W},\tilde{X}%
;R)\rightarrow 0
\end{equation*}%
be the chain of relative complexes (for details, see Section \ref{hml}).
Since 
\begin{equation*}
H_{0}(X;R)\cong R/\langle \alpha (g)x-x\mid g\in \pi ,x\in R\rangle
\rightarrow H_{0}(W;R)\cong R/\langle gx-x\mid g\in G,x\in R\rangle
\end{equation*}%
is surjective, $H_{0}(W,X;R)=0.$ Therefore, the following sequences 
\begin{eqnarray*}
0 &\rightarrow &H_{2}(W,X;R)\rightarrow C_{2}(\tilde{W},\tilde{X}%
;R)\rightarrow \func{Im}i\rightarrow 0; \\
0 &\rightarrow &\func{Im}i\rightarrow \ker j\rightarrow
H_{1}(W,X;R)\rightarrow 0; \\
0 &\rightarrow &\ker j\rightarrow C_{1}(\tilde{W},\tilde{X};R)\rightarrow
C_{0}(\tilde{W},\tilde{X};R)\rightarrow 0
\end{eqnarray*}%
are exact. When $R$ is a principal ideal domain, $H_{2}(W,X;R)$ is always a
free $R$-module, viewed as a submodule of the free module $C_{2}(\tilde{W},%
\tilde{X};R)$. When $R$ is not a principal ideal domain, $%
H_{1}(W,X;R)=H_{1}(G,\pi ;R)$ is a stably free $R$-module by assumption.
According to the exact sequences above, $H_{2}(W,X;R)$ is also stably free
as an $R$-module. By wedging $W$ with some $2$-spheres, which does not
change the fundamental group $G,$ we can further assume that $H_{2}(W,X;R)$
is a free $R$-module.

Since $H_{2}(\tilde{W})\tbigotimes_{\mathbb{Z}[G]}R\rightarrow H_{2}(W,X;R)$
is surjective and $H_{2}(W,X;R)$ is a free $R$-module, by the definition of $%
G$-dense rings we have a set $S$ of elements in $\pi _{2}(W)=H_{2}(\tilde{W}%
) $ whose image forms a basis for $H_{2}(W,X;R).$ Then there are maps $%
b_{\lambda }:S_{\lambda }^{2}\rightarrow W$ with $\lambda \in S$ such that
for all $q\geq 2,$ the composition of maps%
\begin{equation*}
H_{q}(\vee _{\lambda \in S}S_{\lambda }^{2};R)\rightarrow
H_{q}(W;R)\rightarrow H_{q}(W,X;R)
\end{equation*}%
is an isomorphism$.$ For each such $\lambda $, attach a 3-cell $%
(D^{3},S^{2}) $ to $W$ along $b_{\lambda }$. Let $Y$ denote the resulting
space. We see that the diagram 
\begin{equation*}
\begin{array}{ccc}
\vee _{\lambda }S^{2} & \longrightarrow & W \\ 
\downarrow & \ulcorner _{\cdot } & \downarrow \\ 
\vee _{\lambda }D^{3} & \longrightarrow & Y%
\end{array}%
\end{equation*}%
is a pushout diagram. By the van Kampen theorem, the fundamental group of $Y$
is still $G.$ Denoting by $H_{\ast }(-)$ the homology groups $H_{\ast
}(-;R), $ we have the following commutative diagram:%
\begin{equation*}
\begin{array}{ccccccc}
\cdots \rightarrow H_{3}(\vee D^{3},\vee S^{2}) & \rightarrow & H_{2}(\vee
S^{2},\mathrm{pt}) & \rightarrow & H_{2}(\vee D^{3},\mathrm{pt}) & 
\rightarrow & H_{2}(\vee D^{3},\vee S^{2}) \\ 
\downarrow &  & \downarrow &  & \downarrow &  & \downarrow \\ 
\cdots \rightarrow H_{3}(Y,W) & \rightarrow & H_{2}(W,X) & \rightarrow & 
H_{2}(Y,X) & \rightarrow & H_{2}(Y,W).%
\end{array}%
\end{equation*}%
By a five lemma argument, for any $q\geq 2$ the relative homology group $%
H_{q}(Y,X;R)=0,$ which shows that $H_{q}(X;R)\cong H_{q}(Y;R).$
\end{proof}

\section{Applications}

\subsection{Quillen's plus-construction}

In this subsection, we show that Quillen's plus-construction is a special
case of Theorem \ref{th1}. For this, recall that a fibration $F\rightarrow E%
\overset{p}{\rightarrow }B$ with connected fiber $F,$ total space $E$ and
base space $B$ is said to be \textit{quasi-nilpotent }if the action of $\pi
_{1}(B)$ on $H_{\ast }(F;\mathbb{Z})$ is nilpotent. The following result was
proved in \cite{berr}.

\begin{lemma}[4.3(xii) in \protect\cite{berr}]
\label{lem}The following properties of a map $f:X\rightarrow Y$ are
equivalent.

(1) $f$ is acyclic, \textsl{i.e.} $\tilde{H}_{\ast }(F_{f};\mathbb{Z})=0$
for the homotopy fiber $F_{f}$ of $f.$

(2) $f$ is quasi-nilpotent and $H_{\ast }(f)$ is an isomorphism.
\end{lemma}

The following result shows Quillen's plus-construction is a special case of
Theorem \ref{th1} when $R=\mathbb{Z}$ and $\ker \alpha $ is perfect.
According to Lemma \ref{lm2}, $\mathbb{Z}$ is a $G$-dense ring.

\begin{proposition}
\label{prop1}Let $R=\mathbb{Z}$. Among all spaces $Y$ formed from $X$ by
adjoining cells of dimension at most 3 and satisfying (1) in Theorem \ref%
{th1}, it is possible to choose one such that the map $g:X\rightarrow Y$ is
Quillen's plus-construction with respect to $\ker (\alpha )$ if and only if $%
\alpha :\pi \rightarrow G$ is surjective and $\mathrm{ker}(\alpha )$ is
perfect.
\end{proposition}

\begin{proof}
Assume $R=\mathbb{Z}$ and the map $g:X\rightarrow Y$ in Theorem \ref{th1} is
Quillen's plus-construction. Then $g$ is acyclic, \textsl{i.e.} $\tilde{H}%
_{\ast }(F_{g};\mathbb{Z})=0$ for the homotopy fiber $F_{g}$ of $g.$\ This
implies that the homotopy fiber $F_{g}$ is an acyclic space. Since $%
H_{1}(F_{g};\mathbb{Z})=0$, we have that $\pi _{1}(F_{g})$ is perfect. By
the fact that $F_{g}\rightarrow X\rightarrow Y$ is a fiber sequence and $%
\tilde{H}_{0}(F_{g})=0$, then $\alpha $ is surjective and $\mathrm{ker}%
(\alpha )=\mathrm{im}[\pi _{1}(F_{g})\rightarrow \pi ]$ is perfect.

Conversely, when $\alpha $ is surjective and $\mathrm{ker}(\alpha )$ is
perfect, let $\bar{X}$ be the covering space of $X$ with fundamental group $%
\pi _{1}(\bar{X})=\mathrm{ker}(\alpha )$. Since $\mathrm{ker}(\alpha )_{%
\mathrm{ab}}=0$, there is a simply connected space $\overline{Y}$ and a map $%
\bar{g}:\bar{X}\rightarrow \bar{Y}$ such that (1) in Theorem \ref{th1} holds
for any $q\geq 1$ with $R=\mathbb{Z}$ and $G$ the trivial group. Since $%
\overline{Y}$ is simply connected, the map $\bar{g}:\bar{X}\rightarrow \bar{Y%
}$ is quasi-nilpotent and thus acyclic by Lemma \ref{lem}. Let $Y$ be the
pushout of the diagram%
\begin{equation*}
\begin{array}{ccc}
\bar{X} & \longrightarrow & \overline{Y} \\ 
\downarrow & \ulcorner _{\cdot } & \downarrow \\ 
X & \longrightarrow & Y.%
\end{array}%
\end{equation*}%
By the van Kampen Theorem, $\pi _{1}(Y)=\pi _{1}(X)/\pi _{1}(\widetilde{X}%
)=\pi /\mathrm{ker}(\alpha )=G.$ According to (4.20) in \cite{berr}, the map 
$X\rightarrow Y$ is still an acyclic map.
\end{proof}

\subsection{Bousfield's integral localization}

In this subsection, we show that the Bousfield's integral localization is a
special case of Theorem \ref{th1} when $R=\mathbb{Z}$ and $\alpha $ is
surjective, as presented in the following lemma. Again note that $\mathbb{Z}$
is a $G$-dense ring by Lemma \ref{lm2}.

\begin{corollary}[\protect\cite{RS}]
\label{cor}Let $X$ be a CW complex with fundamental group $\pi $ and $N$ a
relatively perfect normal subgroup of $\pi ,$ \textsl{i.e.} $[\pi ,N]=N.$
Then there is a CW complex $Y$ obtained from $X$ by adding $2$-cells and $3$%
-cells such that $\pi _{1}(Y)=\pi /N$ and for any $q\geq 0$ we have%
\begin{equation*}
H_{q}(X;\mathbb{Z})\cong H_{q}(Y;\mathbb{Z}).
\end{equation*}%
Conversely, for some CW complex $Y$ let $f:X\rightarrow Y$ be an integral
homology equivalence of spaces that induces an epimorphism on the
fundamental groups. Then $\ker [f_{\ast }:\pi _{1}(X)\rightarrow \pi
_{1}(Y)] $ is relatively perfect in $\pi _{1}(X).$
\end{corollary}

\begin{proof}
By \cite{h}, there is a long exact sequence%
\begin{equation}
H_{2}(\pi ;\mathbb{Z})\rightarrow H_{2}(\pi /N;\mathbb{Z})\rightarrow N/[\pi
,N]\rightarrow H_{1}(\pi ;\mathbb{Z})\rightarrow H_{1}(\pi /N;\mathbb{Z}%
)\rightarrow 0.  \tag{2}  \label{f2}
\end{equation}%
When $N=[\pi ,N]$, we have that the map $H_{2}(\pi ;\mathbb{Z})\rightarrow
H_{2}(\pi /N)$ is surjective and $H_{1}(\pi ;\mathbb{Z})\rightarrow
H_{1}(\pi /N)$ is an isomorphism. According to Theorem \ref{th1} with $R=%
\mathbb{Z}$, there exists a CW complex $Y$ and a map $g:X\rightarrow Y$ such
that for any $q\geq 0$ we have $H_{q}(X;\mathbb{Z})\cong H_{q}(Y;\mathbb{Z}%
). $ By the proof of Theorem \ref{th1}, $Y$ is obtained from $X$ by
attaching $2 $-cells and $3$-cells. Now assume that $f:X\rightarrow Y$ is an
integral homology equivalence that induces an epimorphism on the fundamental
groups. Then 
\begin{equation*}
H_{1}(X;\mathbb{Z})=H_{1}(\pi ;\mathbb{Z})\cong H_{1}(\pi _{1}(Y);\mathbb{Z}%
)=H_{1}(Y;\mathbb{Z}).
\end{equation*}%
There is a commutative diagram%
\begin{equation*}
\begin{array}{ccc}
H_{2}(X;\mathbb{Z}) & \rightarrow & H_{2}(\pi ;\mathbb{Z}) \\ 
\downarrow &  & \downarrow \\ 
H_{2}(Y;\mathbb{Z}) & \twoheadrightarrow & H_{2}(\pi _{1}(Y);\mathbb{Z}),%
\end{array}%
\end{equation*}%
where the left vertical map is an isomorphism. This shows that the right
vertical map is an epimorphism. According to the same long exact sequence (%
\textsl{\ref{f2}}) above, $\ker f=[\pi ,\ker f],$ \textsl{i.e.}, $\ker f$ is
relatively perfect.
\end{proof}

\subsection{Moore spaces}

In this subsection, we give an application of Theorem \ref{th1} to the
existence of Moore spaces. Recall the definition of Moore spaces from the
Introduction. For $n\geq 2,$ such a space always exists. For $n=1,$ we get
the following result, which was first proved by Varadarajan in \cite{ann}.

When $\pi =1$ in Theorem \ref{th1}, we have the following condition for
existence of $M(G,1)$.

\begin{corollary}[\protect\cite{ann}]
\label{cor2}There exists a Moore space $M(G,1)$ if and only if $H_{2}(G;%
\mathbb{Z})=0.$
\end{corollary}

\begin{proof}
According to Lemma \ref{lm2}, $\mathbb{Z}$ is a $G$-dense ring. If $H_{2}(G;%
\mathbb{Z})=0,$ let $X=\mathrm{pt}$ and $\alpha :1\rightarrow G$ be the
trivial group homomorphism in Theorem \ref{th1}. It is clear that $\alpha $
induces a surjection on $2$-dimensional homology groups and an injection on $%
1$-dimensional homology groups. Then there exists a space $Y$ with $\pi
_{1}(Y)=G$ and $H_{i}(Y;\mathbb{Z})=0$ for any $i\geq 2.$ This makes $Y$ an $%
M(G,1).$

Conversely, let $Y=M(G,1)$ be a Moore space. According to the Hopf exact
sequence (\textsl{cf}. Lemma \ref{lem1})%
\begin{equation*}
\pi _{2}(Y)\rightarrow H_{2}(Y;\mathbb{Z})\rightarrow H_{2}(\pi _{1}(Y);%
\mathbb{Z})\rightarrow 0,
\end{equation*}%
we now have $H_{2}(G;\mathbb{Z})=H_{2}(\pi _{1}(Y);\mathbb{Z})=0.$
\end{proof}

\begin{remark}
It has been noted (\textsl{cf}. \cite{mp}) that the plus-construction shares
some common features with the construction of Moore spaces in \cite{ann}. We
have actually shown in Proposition \ref{prop1} and Corollary \ref{cor2} that
both the plus-construction and the existence of $M(G,1)$ are just two
extreme cases of Theorem \ref{th1} with $\mathbb{Z}$ as coefficients.
\end{remark}

Similarly, we can consider Moore spaces with coefficients. Suppose $G$ is a
group and $R$ is a $\mathbb{Z}[G]$-module. Let $M=M(G,1;R)$ be the Moore
space with coefficients $R,$ \textsl{i.e.} $\pi _{0}(M)=0,\pi _{1}(M)=G$ and 
$H_{i}(M;R)=0$ for $i>1.$ The following new result characterizes when there
exists a Moore space with coefficients $R.$

\begin{proposition}
Let $R$ be a $G$-dense ring. Suppose that $R$ is a principal ideal domain or
that the relative homology group $H_{1}(G,1;R)$ is a stably free $R$-module.
Then there exists a Moore space $M(G,1;R)$ if and only if $H_{2}(G;R)=0.$
\end{proposition}

\begin{proof}
The proof is similar to the proof of Corollary \ref{cor2}.
\end{proof}

\subsection{Partial $k$-completion of Bousfield and Kan}

Let $k=\mathbb{Z}/p$ be the constant ring for any prime $p$ or $k\subseteq 
\mathbb{Q}$ a subring of the rationals. A group $G$ is called $k$-perfect,
if $H_{1}(G;k)=k\bigotimes_{\mathbb{Z}}G_{ab}=0.$ The following result was
first obtained by Bousfield and Kan \cite{Bk}; it is a special case of
Theorem \ref{th1} when $\alpha $ is surjective.

\begin{corollary}[Prop. 6.3 in \protect\cite{Bk}, p. 219]
\label{cmp}Assume $k=\mathbb{Z}/p$ is the constant ring for any prime $p$ or 
$k\subseteq \mathbb{Q}$ a subring of rationals. Let $X$ be a CW complex with
fundamental group $\pi $ with $P$ its maximal $k$-perfect subgroup. Then
there exists a CW complex $Y$ and a map $g:X\rightarrow Y$ such that $\pi
_{1}(Y)=\pi /P$ and for any $q\geq 0$ we have 
\begin{equation*}
H_{q}(X;k[\pi /P])\cong H_{q}(Y;k[\pi /P]).
\end{equation*}
\end{corollary}

\begin{proof}
By \cite{h}, there is a long exact sequence%
\begin{eqnarray*}
H_{2}(\pi ;k[\pi /P]) &\rightarrow &H_{2}(\pi /P;k[\pi /P])\rightarrow k[\pi
/P]\tbigotimes\nolimits_{\mathbb{Z}[\pi /P]}P_{\mathrm{ab}} \\
&\rightarrow &H_{1}(\pi ;k[\pi /P])\rightarrow H_{1}(\pi /P;k[\pi
/P])\rightarrow 0.
\end{eqnarray*}%
When $k[\pi /P]\tbigotimes_{\mathbb{Z}[\pi /P]}P_{\mathrm{ab}}\cong
k\bigotimes_{\mathbb{Z}}P_{\mathrm{ab}}=0$, we can see 
\begin{equation*}
H_{2}(\pi ;k[\pi /P])\rightarrow H_{2}(\pi /P;k[\pi /P])
\end{equation*}%
is surjective and $H_{1}(\pi ;k[\pi /P])\rightarrow H_{1}(\pi /P;k[\pi /P])$
is an isomorphism. According to the long exact sequence of homology groups%
\begin{eqnarray*}
\cdots &\rightarrow &H_{1}(\pi ;k[\pi /P])\rightarrow H_{1}(\pi /P;k[\pi
/P])\rightarrow H_{1}(\pi /P,\pi ;k[\pi /P]) \\
&\rightarrow &H_{0}(\pi ;k[\pi /P])\rightarrow H_{0}(\pi /P;k[\pi
/P])\rightarrow 0,
\end{eqnarray*}%
we have $H_{1}(\pi /P,\pi ;k[\pi /P])=0.$ Therefore, applying Theorem \ref%
{th1} with $G=\pi /P$ and the $G$-dense ring $R=k[\pi /P],$ we get a CW
complex $Y$ and a map $g:X\rightarrow Y$ such that for any $q\geq 0$ there
is an isomorphism 
\begin{equation*}
H_{q}(X;k[\pi /P])\cong H_{q}(Y;k[\pi /P]).
\end{equation*}
\end{proof}

\subsection{Zero-in-the-spectrum conjecture\label{hml}}

In this subsection, we give an application of Theorem \ref{th1} to the
zero-in-the-spectrum conjecture. This conjecture is stated in $L^{2}$%
-homology, whose definition and basic properties are presented as
follows.\medskip

\noindent $L^{2}$\textbf{-homology.}

Let $G$ be a group and $l^{2}(G)$ be the Hilbert space spanned by $G$ (%
\textsl{cf}. Section $\ref{s1}$). By definition, the reduced group $C^{\ast
} $-algebra $C_{r}^{\ast }(G)$ ( resp. $C_{\mathbb{R}}^{\ast }(G)$) is the
completion of $\mathbb{C}[G]$ (resp. $\mathbb{R}[G]$) in $B(l^{2}(G))$ with
respect to the operator norm. Suppose $Y$ is a $G$-CW complex and $C_{\ast
}(Y)$ the cellular chain complex of $Y.$ The $L^{2}$ $n$-th homology group
with coefficients $C_{r}^{\ast }(G)$ of $Y$ is defined as $H_{n}(C_{\ast
}(Y)\bigotimes_{\mathbb{Z}[G]}C_{r}^{\ast }(G)).$ When $Y$ is a CW complex
with fundamental group $G,$ the universal covering space $\tilde{Y}$ is a $G$%
-CW complex with the group action as deck transformation. The $L^{2}$%
-homology groups with coefficients $C_{r}^{\ast }(G)$ of $Y$ are defined as
the homology groups of the complex $C_{\ast }(\tilde{Y})\bigotimes_{\mathbb{Z%
}[G]}$ $C_{r}^{\ast }(G),$ where $C_{\ast }(\tilde{Y})$ is the cellular
chain complex of $\tilde{Y}$.

Let $f:X\rightarrow Y$ be a cellular map of CW complexes and assume $\pi
_{1}(Y)=G$. We view $C_{i}(\tilde{X})$ as a free right $\mathbb{Z}[\pi
_{1}(X)]$-module spanned by the $i$-th dimensional cells in $X$ and define
the $L^{2}$-homology groups of $X$ with coefficients $C_{r}^{\ast }(G)$ as 
\begin{eqnarray*}
H_{\ast }(X;C_{r}^{\ast }(G)) &=&H_{\ast }((C_{\ast }(\tilde{X}%
)\tbigotimes\nolimits_{\mathbb{Z}[\pi _{1}(X)]}\mathbb{Z}[G])\tbigotimes%
\nolimits_{\mathbb{Z}[G]}C_{r}^{\ast }(G)) \\
&=&H_{\ast }(C_{\ast }(\tilde{X})\tbigotimes\nolimits_{\mathbb{Z}[\pi
_{1}(X)]}C_{r}^{\ast }(G)).
\end{eqnarray*}%
Suppose $f_{1}:\pi _{1}(X)\rightarrow \pi _{1}(Y)$ is the group homomorphism
induced by $f.$ Then there is a well-defined chain map%
\begin{equation*}
f_{\ast }:C_{\ast }(\tilde{X})\tbigotimes\nolimits_{\mathbb{Z}[\pi _{1}(X)]}%
\mathbb{Z}[G]\rightarrow C_{\ast }(\tilde{Y}),\text{ }xh\tbigotimes g\mapsto
f(x)f_{1}(h)g
\end{equation*}%
for any cell $x\in X,h\in \pi _{1}(X)$ and $g\in G.$ Taking the mapping
cylinder if necessary, we can assume $f_{\ast }$ is an inclusion of free $%
\mathbb{Z}[G]$-modules. We can define the relative $L^{2}$-homology as the
homology of the cokernel chain complex $C_{\ast }(\tilde{Y},\tilde{X}%
;C_{r}^{\ast }(G))$ of the chain map $f_{\ast }$ in the usual sense to get a
long exact sequence%
\begin{eqnarray*}
\cdots &\rightarrow &H_{n}(X;C_{r}^{\ast }(G))\rightarrow
H_{n}(Y;C_{r}^{\ast }(G))\rightarrow H_{n}(Y,X;C_{r}^{\ast }(G)) \\
&\rightarrow &H_{n-1}(X;C_{r}^{\ast }(G))\rightarrow H_{n-1}(Y;C_{r}^{\ast
}(G))\cdots \rightarrow H_{0}(Y,X;C_{r}^{\ast }(G))\rightarrow 0.
\end{eqnarray*}%
When we consider the group von Neumann algebra $\mathcal{N}G$ or the Hilbert
space $l^{2}(G)$ instead of the reduced group $C^{\ast }$-algebra $%
C_{r}^{\ast }(G),$ the $L^{2}$-homology with coefficients $\mathcal{N}G$ or $%
l^{2}(G)$ can be defined in a similar way. For more details on $L^{2}$%
-homology, we refer the reader to the book of L\"{u}ck \cite{lu} and the
article \cite{hig}.

\begin{proposition}
\label{prop2}Let $X$ be a CW complex with fundamental group $G.$ Assume $%
R^{\prime }\subseteq R$ are two $\mathbb{Z}[G]$-modules with $R^{\prime }$ a
unital ring such that $R$ is an $R^{\prime }$-module and $R/R^{\prime }$ is
a flat $R^{\prime }$-module. Then for any nonnegative integer $n,$ the
following are equivalent:

\begin{enumerate}
\item[(i)] For any $0\leq i\leq n,$ the homology group $H_{i}(X;R^{\prime
})=0.$

\item[(ii)] For any $0\leq i\leq n,$ the homology group $H_{i}(X;R)=0$.
\end{enumerate}

In particular, the homology groups $H_{i}(X;C_{\mathbb{R}}^{\ast }(G))$ are
zero in degrees $0$ through $n$ if and only if the homology groups $%
H_{i}(X;C_{r}^{\ast }(G))$ are zero in degrees $0$ through $n.$
\end{proposition}

\begin{proof}
Suppose that for any $0\leq i\leq n$, we have $H_{i}(X;R^{\prime })=0.$ Let $%
(C_{\ast }(\tilde{X}))$ be the chain complex of the universal covering space 
$\tilde{X},$ which is viewed as a chain complex of free $\mathbb{Z}[G]$%
-modules. Thus $C_{\ast }(\tilde{X})\bigotimes_{\mathbb{Z}[G]}R^{\prime }$
is a chain complex of free $R^{\prime }$-modules. According to the K\"{u}%
nneth spectral sequence (\textsl{cf}. page 143 in \cite{we}), we have 
\begin{equation*}
\mathrm{Tor}_{p}^{R^{\prime }}(H_{q}(X;R^{\prime }),R)\Rightarrow
H_{p+q}(X;R),
\end{equation*}%
which shows that for any $0\leq i\leq n,$ the homology group $H_{i}(X;R)=0.$

Conversely, assume that any for $0\leq i\leq n,$ we have $H_{i}(X;R)=0.$
Since $\mathrm{Tor}_{0}^{R^{\prime }}(M,N)=M\bigotimes_{R^{\prime }}N$ for
any $R^{\prime }$ -modules $M$ and $N,$ 
\begin{equation*}
H_{0}(X;R^{\prime })\tbigotimes\nolimits_{R^{\prime }}R=H_{0}(X;R)=0
\end{equation*}%
by the K\"{u}nneth spectral sequence. According to the long exact sequence%
\begin{eqnarray*}
\cdots &\rightarrow &\mathrm{Tor}_{1}^{R^{\prime }}(H_{0}(X;R^{\prime
}),R/R^{\prime })\rightarrow H_{0}(X,R^{\prime
})\tbigotimes\nolimits_{R^{\prime }}R^{\prime }=H_{0}(X;R^{\prime }) \\
&\rightarrow &H_{0}(X,R^{\prime })\tbigotimes\nolimits_{R^{\prime
}}R=0\rightarrow H_{0}(X,R^{\prime })\tbigotimes\nolimits_{R^{\prime
}}R/R^{\prime }\rightarrow \cdots ,
\end{eqnarray*}%
there is a surjection 
\begin{equation*}
\mathrm{Tor}_{1}^{R^{\prime }}(H_{0}(X;R^{\prime }),R/R^{\prime
})\twoheadrightarrow H_{0}(X;R^{\prime }).
\end{equation*}%
Since $R/R^{\prime }$ is a flat $R^{\prime }$-module by assumption, $%
H_{0}(X;R^{\prime })=0.$ By the same spectral sequence, we can get 
\begin{equation*}
\mathrm{Tor}_{0}^{R^{\prime }}(H_{1}(X;R^{\prime }),R)=H_{1}(X;R)=0.
\end{equation*}%
Similar arguments prove inductively that for any $0\leq i\leq n,\ $the
homology group $H_{i}(X;R^{\prime })=0$. Clearly, $C_{r}^{\ast }(G)\cong C_{%
\mathbb{R}}^{\ast }(G)\bigoplus iC_{\mathbb{R}}^{\ast }(G)$ if we consider
the real and imaginary parts. The fact that $C_{r}^{\ast }(G)/C_{\mathbb{R}%
}^{\ast }(G)$ is a free $C_{\mathbb{R}}^{\ast }(G)$-module proves the last
part of the proposition.\medskip
\end{proof}

\noindent \textbf{Zero-in-the-spectrum conjecture.}

In the notation of the Introduction, \emph{zero not belonging to the
spectrum }of\emph{\ }$\Delta =\Delta _{\ast }$\emph{\ }can also be expressed
as the vanishing of $H_{\ast }(M;C_{r}^{\ast }(\pi _{1}(M))).$ For more
details, we refer the reader to \cite{lu}.

\begin{conjecture}
Let $M$ be a finite aspherical $CW$-complex (or weakly a closed, connected,
oriented and aspherical Riemannian manifold) with fundamental group $\pi .$
Then for some $i\geq 0,$ $H_{i}(X;C_{r}^{\ast }(\pi ))\neq 0.$
\end{conjecture}

If the condition that $X$ is aspherical is dropped, the following corollary,
which is a special case of Theorem \ref{th1} when $R=C_{\mathbb{R}}^{\ast
}(G)\ $and $\pi =1$, shows the above conjecture is not true. This result is
a generalization of the results obtained by Farber-Weinberger \cite{FW} and
Higson-Roe-Schick \cite{hig}. Here we do not assume $G$ is finitely
presented. Recall from Lemma \ref{lm2} that $C_{\mathbb{R}}^{\ast }(G)$ is a 
$G$-dense ring.

\begin{corollary}[\protect\cite{hig}]
\label{zero}For a group $G$ with 
\begin{equation*}
H_{0}(G;C_{r}^{\ast }(G))=H_{1}(G;C_{r}^{\ast }(G))=H_{2}(G;C_{r}^{\ast
}(G))=0,
\end{equation*}
there is a CW complex $Y$ such that $\pi _{1}(Y)=G$ and for each integer $%
n\geq 0,$ the homology group $H_{n}(Y;C_{r}^{\ast }(G))=0.$ If $G$ is
finitely presented, $Y$ can be chosen to be a finite CW complex.
\end{corollary}

\begin{proof}
According to Proposition \ref{prop2}, 
\begin{equation*}
H_{0}(G;C_{\mathbb{R}}^{\ast }(G))=H_{1}(G;C_{\mathbb{R}}^{\ast
}(G))=H_{2}(G;C_{\mathbb{R}}^{\ast }(G))=0.
\end{equation*}
Let $\alpha :\pi =1\rightarrow G,$ $R=C_{\mathbb{R}}^{\ast }(G)$ and $X=%
\mathrm{pt}$ in Theorem \ref{th1}. By the long exact sequence of homology
groups%
\begin{eqnarray*}
\cdots &\rightarrow &H_{1}(1;C_{\mathbb{R}}^{\ast }(G))\rightarrow
H_{1}(G;C_{\mathbb{R}}^{\ast }(G))\rightarrow H_{1}(G,1;k[\pi /P]) \\
&\rightarrow &H_{0}(1;C_{\mathbb{R}}^{\ast }(G))\rightarrow H_{0}(G;C_{%
\mathbb{R}}^{\ast }(G))\rightarrow 0,
\end{eqnarray*}%
we have 
\begin{equation*}
H_{1}(G,\pi ;C_{\mathbb{R}}^{\ast }(G))=H_{0}(1;C_{\mathbb{R}}^{\ast
}(G))=C_{\mathbb{R}}^{\ast }(G),
\end{equation*}%
which is a free $C_{\mathbb{R}}^{\ast }(G)$-module. Therefore, there exists
such $Y$ by Theorem \ref{th1} and Proposition \ref{prop2} such that $\pi
_{1}(Y)=G$ and for all $n\geq 0$ we have $H_{n}(Y;C_{r}^{\ast }(G))=0$. When 
$G$ is finitely presented, the proof of Theorem \ref{th1} shows the number
of cells added to $X$ is finite.
\end{proof}

The groups satisfying the condition of Corollary \ref{zero} have to be
non-amenable (for more details, see \cite{Br} and Rosenberg's note in \cite%
{ros}).

\bigskip

\noindent \textbf{Acknowledgements}

The author is grateful to his advisor Professor A. J. Berrick for
introducing this subject to him and for many helpful discussions. Helpful
input of the referee is also acknowledged.

\bigskip

Department of Mathematics, National University of Singapore, Kent Ridge
119076, Singapore.

E-mail: yeshengkui@nus.edu.sg

\end{document}